\documentclass[10pt]{article}
\usepackage[margin=1in]{geometry}

\usepackage{amsmath, amssymb, amsthm}
\usepackage[mathic]{mathtools}
\usepackage{algorithmic}
\usepackage{algorithm}
\usepackage[T1]{fontenc}
\usepackage{bbm}
\usepackage{float}
\usepackage{subcaption}
\usepackage{booktabs}
\usepackage{libertine}
\usepackage{authblk}
\usepackage{algorithmic}
\usepackage{algorithm}
\usepackage{optidef}
\usepackage{relsize}
\usepackage[toc,page]{appendix}
\usepackage{listings}
	
\usepackage{lineno}

\usepackage{thm-restate}
\usepackage[libertine]{newtxmath}
\usepackage[scaled=0.96]{zi4} 

\usepackage{a4wide, dsfont, microtype, xcolor, paralist}

\usepackage[ocgcolorlinks]{hyperref} 
\colorlet{DarkRed}{red!50!black}
\colorlet{DarkGreen}{green!50!black}
\colorlet{DarkBlue}{blue!50!black}

\theoremstyle{thmstyleone}%
\newtheorem{theorem}{Theorem}[section]
\newtheorem{definition}{Definition}[section]%
\newtheorem{lemma}{Lemma}[section]
\newtheorem{conjecture}{Conjecture}[section]
\newtheorem{proposition}[theorem]{Proposition}%

\title{Rectilinear crossing number of the double circular complete bipartite graph}
\author{Mohsen Nafar}
\affil{Bielefeld University \\
mohsen.nafar@uni-bielefeld.de}
\date{}

\begin{document}
\maketitle

\begin{abstract}
In this work, we study a mathematically rigorous metric of a graph visualization quality under conditions that relate to visualizing a bipartite graph. Namely we study rectilinear crossing number in a special arrangement of the complete bipartite graph (\(K_{m,n}\), i.e., \textit{Double circular arrangement} (\(K_{m,n}^{\mathbf{DC}}\)) where the two parts are placed on two concentric circles. For this purpose, we introduce a combinatorial formulation to count the number of crossings. We prove a proposition about the rectilinear crossing number of \(K_{m,n}^{\mathbf{DC}}\). Then, we introduce a geometric optimization problem whose solution gives the optimum radii ratio in the case that the number of crossings for \(K_{m,n}^{\mathbf{DC}}\) is minimized. Later on, we study the magnitude of change in the number of crossings upon change in the radii of the circles. In this part, we present and prove a lemma on bounding the changes in the number of crossings of \(K_{m,n}^{\mathbf{DC}}\) that is followed by a theorem on asymptotics of the bounds.
\end{abstract}


\textbf{keywords:} Rectilinear crossing number of \(K_{n,m}\), Crossing number bounds, Complete bipartite graph crossing changes



\section{Introduction}\label{sec1}

The problem of the crossing number in graphs has been studied for decades \cite{zarankiewicz1955problem, schaefer2012graph, pach1996applications, garey1983crossing}. One area of application for which the study of crossing number is useful is graph visualization where lowering of the number of the edge crossings can improve the readability of the resulting layout. 

In a separate work we have presented an algorithm for hypergraph visualization in which we visualize the bipartite representation of a co-authorship network. In the bipartite representation of the network one part of the vertices corresponds to authors and the other part represents the articles so as to show the relations in the network such as relations between authors with each other and relations between authors and articles. In our algorithm, one of our objectives is to minimize the number of crossings to improve the readability of the layout \cite{our-work}. In our visualization algorithm, in the first stage of the algorithm, we place the vertices of the bipartite representation of a hypergraph on two concentric circles (i.e. one part of the bipartite graph vertices on one circle and the other part on the other circle). Therefore, it will be a good idea to try to start the first stage of an algorithm that minimizes the number of crossings as one of its objectives, by positioning the vertices on positions such that the bipartite graph has the minimum number of crossings. 

We study the complete bipartite graph since amongst all bipartite graphs it possesses the fewest number of parameters and is the most symmetric one. Moreover, obtaining an upper bound for the number of crossings of a complete bipartite graph also works for bipartite graphs in general. Furthermore, a lower bound on the number of crossings for complete bipartite graph works as the worst case for bipartite graphs in general. And at the same time the crossing number of a complete bipartite graph is the worst case in terms of the number of crossings because it contains all the edges. 

The points mentioned above and the algorithmic challenges of hypergraph visualization motivated the study of the rectilinear crossing number of the complete bipartite graph in a special arrangement of the vertices on plane. In this specific embedding of the complete bipartite graph on plane the vertices of the two parts are placed on two concentric circles. We call this arrangement of the complete bipartite graph on plane, \textit{double circular} and denote it by \(K_{m,n}^{\mathbf{DC}}\) (see Figure~\ref{double circular graphs}). In this paper, first, we prove the following proposition:

\begin{figure}[H]
    \centering
    \includegraphics[width=8cm]{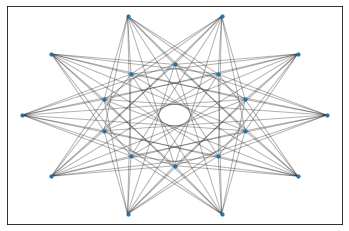}
    \caption{Double circular arrangement of \(K_{10,10}\): \(K_{10,10}^{\mathbf{DC}}\)}
    \label{double circular graphs}
\end{figure}

\begin{proposition}\label{Prop1}
    Let \(m \ge n\), the rectilinear crossing number of the double circular arrangement of \(K_{m,n}\) (see ) is

\begin{equation*}
    cr\Big(K_{m,n}^{\mathbf{DC}}\Big) = \binom{m}{2} \cdot \bigg \lfloor{\frac{n}{2}}\bigg\rfloor \cdot \bigg \lfloor{\frac{n-1}{2}}\bigg \rfloor.
\end{equation*}
\end{proposition}

Then, we define a geometric optimization problem whose solution results in a threshold for the ratio between the radii of the circles up to which the number of the rectilinear crossings of \(K_{m,n}^\mathbf{DC}\) remains equal to its rectilinear crossing number, which we stated in the previous proposition. 

Later, we prove the following lemma on bounding the changes in the number of the rectilinear crossings of \(K_{m,n}^\mathbf{DC}\) when the ratio between the radii of the two circles exceeds the threshold. We denote it by \(N_{cr}\Big(K_{m,n}^{\mathbf{DC}}\Big)\).

\begin{lemma}\label{Lemma1}
Let \(R\) and \(r\) be the radii of outer and inner circles in \(K_{m,n}^{\mathbf{DC}}\) arrangement, respectively. Moreover, assume that the parts with \(m\) and \(n\) vertices are placed on the inner and the outer circle, respectively. Let \(cr\Big(K_{m,n}^{\mathbf{DC}}\Big)\) be the rectilinear crossing number for \(K_{m,n}^{\mathbf{DC}}\), then we have the following bounds for \(N_{cr}\Big(K_{m,n}^{\mathbf{DC}}\Big)\):

\begin{equation*}
    \begin{cases}
        N_{cr}\Big(K_{m,n}^{\mathbf{DC}}\Big) - cr\Big(K_{m,n}^{\mathbf{DC}}\Big) \ge m \sum_{j=1}^{\lceil\frac{m}{2}\rceil-1} \Bigg(\Big(\frac{1-(-1)^n}{2}\Big)\beta_j + (\beta_j)^2\Bigg)   
        \text{ }\\\\
        N_{cr}\Big(K_{m,n}^{\mathbf{DC}}\Big) - cr\Big(K_{m,n}^{\mathbf{DC}}\Big) \le 4m  \sum_{j=1}^{\lceil\frac{m}{2}\rceil-1} \Bigg(\Big(\frac{1-(-1)^n}{4}\Big)\beta_j + (\beta_j)^2\Bigg),

    \end{cases}
\end{equation*}

where
    \begin{equation*}\label{beta-theorem}
    \beta_j = \Bigg \lfloor {\frac{r\cdot \sin{\Big(\frac{j\cdot \pi}{m}\Big)}}{2\cdot R\cdot \sin{\big(\frac{\pi}{n}\big)}}} \Bigg \rfloor.
    \end{equation*}
\end{lemma}

Note that \(\beta_j=\Theta\left(\frac{r\cdot\Theta(1)}{R\cdot\Theta(1/n)}\right)=\Theta(n)\) for \(n\to\infty\), \(r/R=\Theta(1)\) and \(j=\Omega(m)\), and thus for such parameters we have 
\[\sum_{j=1}^{\lceil\frac{m}{2}\rceil-1} \Bigg(\Big(\frac{1-(-1)^n}{4}\Big)\beta_j + (\beta_j)^2\Bigg)=\Theta(n^2),\]
which results in the following asymptotic bound for the change in the number of rectilinear crossings:


\begin{theorem}
    \[N_{cr}\Big(K_{m,n}^{\mathbf{DC}}\Big)=\binom{m}{2} \bigg\lfloor\frac{n}{2}\bigg\rfloor  \bigg\lfloor\frac{n-1}{2}\bigg\rfloor+\Theta(n^2m).\]
\end{theorem}

An arrangement similar to this is the plane drawing of the one called \textit{cylindrical drawing} of a bipartite graph. Richter and Thomassen in \cite{richter1997relations} proved that the crossing number of the complete bipartite graph \(K_{n,n}\) in a cylindrical drawing is
\begin{equation*}
    cr_{\circledcirc}(K_{n,n}) = n \binom{n}{3}.
\end{equation*}

The differences between a cylindrical drawing and a double circular arrangement are the followings. First, edges in cylindrical drawing are not restricted to be straight line segments whereas in double circular arrangement edges are straight line segments. Second, edges in cylindrical drawing cannot cut the inner circle. Abrego, Fernandez-Merchant, and Sparks \cite{abrego2020bipartite} proved the two following theorems on the crossing number of the complete bipartite graph in cylindrical drawing, which they denote by \(cr_{\circledcirc}(K_{m,n})\). 

\begin{theorem}
    If \(m \le n\), the bipartite cylindrical crossing number of \(K_{m,n}\) is
\begin{multline}
        cr_{\circledcirc}(K_{m,n}) = \binom{n}{2} \binom{m}{2} + \sum_{1\le i \le j \le m} \Bigg(\bigg \lfloor \frac{n}{m} (j-1) \bigg \rfloor - \bigg \lfloor \frac{n}{m} (i-1) \bigg \rfloor \Bigg)^2  \\ - n \sum_{1\le i \le j \le m} \Bigg(\bigg \lfloor \frac{n}{m} (j-1) \bigg \rfloor - \bigg \lfloor \frac{n}{m} (i-1) \bigg \rfloor \Bigg).
\end{multline}
\end{theorem}

\begin{theorem}
    If \(m\) divides \(n\), the bipartite cylindrical crossing number of \(K_{m,n}\) is
    \begin{equation*}
    cr_{\circledcirc}(K_{m,n}) = \frac{1}{12} n (m-1) (2mn - 3m - n).
\end{equation*}
\end{theorem}

\section{Double circular arrangement of Complete bipartite graph}\label{Arrangements}

In this section we will move gradually from the simplest case of arrangement of the complete bipartite graph on plane towards the Double circular arrangement. From now on we denote the two parts of the vertices of the complete bipartite graph \(K_{m,n}\) with part \(A\) and part \(B\).

\subsection{First arrangement: Double Parallel Lines DPL}\label{DPL-arrangement}

\begin{figure}[H]
    \centering
    \includegraphics[width=6cm]{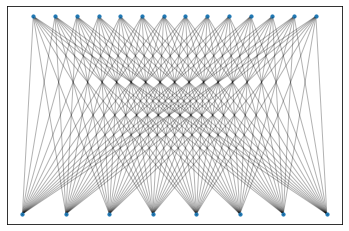}
    \caption{Double parallel lines arrangement of \(K_{14,8}:K_{14,8}^{\mathbf{DPL}}\)}
    \label{fig_DPL}
\end{figure}

Suppose that the two parts of \(K_{m,n}\) are placed on two parallel lines and denote it by \(K_{m,n}^{\mathbf{DPL}}\) (i.e. part \(A\) on one line and part \(B\) on the other, see Figure~\ref{fig_DPL}). We can count the total number of crossings by computing the edge crossings caused by edges between every pair \(\{i,j\}\) of vertices of part \(A\) and all vertices of the part \(B\) which we denote by \(cr_{\{i,j\}}^n\)- without loss of generality assume that the part \(A\) is the part with \(m\) vertices- and then we sum up these numbers for all pairs. Therefore, we count the number of crossings as the following:

for every pair \(\{i,j\}\) that belongs to the part \(A\) we have one \(K_{2,n}^{\mathbf{DPL}}\) then we compute its edge crossings using the following formula:

\begin{equation*}
    cr_{\{i,j\}}^n = cr(K_{2,n}^{\mathbf{DPL}}) = \sum_{k=0}^{n-1} k = \binom{n}{2},
\end{equation*}

since there are \(\binom{m}{2}\) pairs of vertices in the part \(A\), we get the following number of crossings for \(K_{m,n}^{\mathbf{DPL}}\):

\begin{equation}\label{tpl-compute}
    cr\Big(K_{m,n}^{ \mathbf{DPL}}\Big) = \binom{m}{2} \cdot \binom{n}{2}.
\end{equation}

\subsection{Second arrangement: Triple Parallel Lines TPL}\label{TPL-arrangement} 

\begin{figure}[H]
    \centering
    \includegraphics[width=6cm]{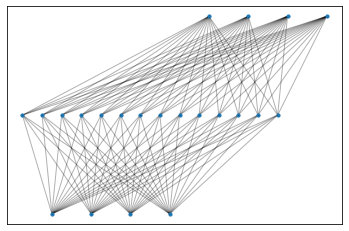}
    \caption{Triple parallel lines arrangement of \(K_{14,8}:K_{14,8}^{(4, 14, 4)}\)}
    \label{fig_tpl}
\end{figure}

Consider three parallel lines \textit{top line}, \textit{center line} and \textit{bottom line}. In this arrangement we put part \(A\) on the center line and divide the part \(B\) into two partitions and place them on the top and bottom line that are on two sides of the center line and are parallel to it (see Figure~\ref{fig_tpl}). We have \(t,b\) and \(c\) vertices on top, bottom and center line, respectively. We denote this arrangement by \(K_{m,n}^{\mathbf{TPL}}\).

Here we introduce the notation \(K_{m,n}^{(t, m, b)}\) that we are going to use in the rest of this section. It means that in \(K_{m,n}^{\mathbf{TPL}}\), the center line contains the part with \(m\) vertices and the \(n\) vertices of the other part are divided into \(t\) and \(b\) vertices and are placed on top and bottom line, respectively. 

Now we count the number of crossings for \(K_{m,n}^{(t, m, b)}\) with the following formula:

for every pair \(\{i,j\}\) of vertices that belong to the center line we have

\begin{equation}\label{one-pair-tpl}
    cr\Big(K_{2,t+b}^{(t, 2, b)}\Big) = cr_{2}^t + cr_{2}^b = \sum_{x=1}^{t-1} x + \sum_{y=1}^{b-1} y = \binom{t}{2} + \binom{b}{2}.
\end{equation}

Finally, since there exist \(\binom{c}{2}\) distinct pairs of vertices in the part on the center line, the number of crossings for \(K_{m,n}^{(t, m, b)}\) is as follows:

\begin{equation}\label{number-cr-tpl}
    cr\Big(K_{c,t+b}^{(t, c, b)}\Big) = \binom{c}{2} \cdot \Bigg(\binom{t}{2} + \binom{b}{2}\Bigg).
\end{equation}

Since the function \(f(x) = x^2 + (z-x)^2\) attains its minimum at the point \(x=\lceil \frac{z}{2}\rceil\) or \(x=\lfloor \frac{z}{2}\rfloor\), it is easy to see that the following inequality is true for all positive integers \(x,y\), where \(x+y=z\):

\begin{equation}\label{binom-min}
    \binom{x}{2} + \binom{y}{2} \ge \binom{\lfloor \frac{z}{2}\rfloor}{2} + \binom{\lceil \frac{z}{2}\rceil}{2}.
\end{equation}

Moreover, it is easy to check that for every pair of positive integers \(x,y\), where \(x \le y\), the following inequality is true:

\begin{equation}\label{m-n-inequality}
    \binom{x}{2} \cdot \bigg \lfloor{\frac{y}{2}}\bigg \rfloor \cdot \bigg \lfloor{\frac{y-1}{2}}\bigg \rfloor \ge \binom{y}{2} \cdot \bigg \lfloor{\frac{x}{2}}\bigg \rfloor \cdot \bigg \lfloor{\frac{x-1}{2}}\bigg \rfloor.     
\end{equation}

Hence, in order to get the minimum number of crossings in this arrangement, without loss of generality assume that \(m \ge n\). By using inequality~\eqref{m-n-inequality} and inequality~\eqref{binom-min} in equation~\eqref{number-cr-tpl}, we get that we should put the part with larger number of vertices on the center line (i.e., \(c = m\)) and the other part of the vertices on the top and bottom lines (i.e. \(\lvert t - b \rvert \le 1\)). Then we get the following for the crossing number of \(K_{m,n}^{\mathbf{TPL}}\), where \(m \ge n\):

\begin{equation*}
    Cr\Big(K_{m,n}^{\mathbf{TPL}}\Big) = \binom{m}{2} \cdot \bigg \lfloor{\frac{n}{2}}\bigg \rfloor \cdot \bigg \lfloor{\frac{n-1}{2}}\bigg \rfloor = Cr\Bigg(K_{m,n}^{\big(\lfloor \frac{n}{2}\rfloor, m, \lceil \frac{n}{2}\rceil \big)}\Bigg).
\end{equation*}

\subsection{Third arrangement: Double Orthogonal Lines DOL}\label{DOL-arrangement}

Double Orthogonal Lines arrangement is the one in which the vertices of \(K_{m,n}\) are placed on two imaginary perpendicular lines. We denote this arrangement by \(K_{m,n}^{\mathbf{DOL}}\). We call these two lines \textit{horizontal} and \textit{vertical} lines and denote them by line \(H\) and line \(V\), respectively. 

Suppose the lines \(H\) and \(V\) contain the part \(A\) with \(m\) and part \(B\) with \(n\) vertices of \(K_{m,n}\), respectively. Let these two perpendicular lines cut each other in a way that \(m\) vertices on the line \(H\) gets divided into two partitions with \(l\) and \(r\) vertices, and \(n\) vertices on the line \(V\) gets divided into two partitions with \(t\) and \(b\) vertices. 

Here we introduce a new notation \(K_{m,n}^{((l,r),(t,b))}\) that we are going to use in the rest of this section. The meaning of it is that considering the arrangement \(K_{m,n}^{\mathbf{DOL}}\), the \(m\) vertices on the horizontal line are divided into two partitions, i.e. \((l,r)\), and the \(n\) vertices on the vertical line are divided into two partitions, i.e. \((t,b)\). 

Then, we can compute the number of crossings for \(K_{m,n}^{((l,r),(t,b))}\) in the following way:

\begin{equation}\label{dol-compute}
    cr\Big(K_{m,n}^{((l,r),(t,b))}\Big) =  cr\Big(K_{l,t}^{\mathbf{DPL}}\Big) +  
    cr\Big(K_{l,b}^{\mathbf{DPL}}\Big) + 
    cr\Big(K_{r,t}^{\mathbf{DPL}}\Big) +  
    cr\Big(K_{r,b}^{\mathbf{DPL}}\Big).
\end{equation}

By using equation~\eqref{tpl-compute} in equation~\eqref{dol-compute} we get the following inequality:

\begin{equation*}
    cr\Big(K_{m,n}^{((l,r),(t,b))}\Big) =  \binom{l}{2}\cdot \binom{t}{2} + 
    \binom{l}{2}\cdot \binom{b}{2} +
    \binom{r}{2} \cdot \binom{t}{2} +
    \binom{r}{2}\cdot \binom{b}{2}.
\end{equation*}

And after we simplify the previous formula we get the following equation:

\begin{equation}\label{dol-compute-last}
    cr\Big(K_{m,n}^{((l,r),(t,b))}\Big) =  
    \Bigg(\binom{l}{2} + \binom{l}{2}\Bigg) \cdot 
    \Bigg(\binom{t}{2} + \binom{b}{2}\Bigg).
\end{equation}

Since the function \(f(x) = x^2 + (n-x)^2\) attains its maximum at point \(x=n\) or \(x=0\), the following inequality is true for all non negative integers \(x,y\), where \(x+y=n\):

\begin{equation}\label{binom-max}
    \binom{x+y}{2} \ge \binom{x}{2} + \binom{y}{2}.
\end{equation}

\begin{figure}[H]
    \centering
    \includegraphics[width=6cm]{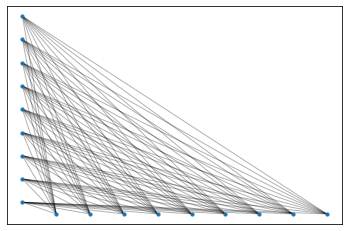}
    \caption{Double orthogonal lines arrangement of \(K_{9,9}:K_{9,9}^{((0,9),(9,0))}\) with maximum number of crossings}
    \label{fig_DOLmax}
\end{figure}

We use inequality~\eqref{binom-max} in equation~\eqref{dol-compute-last} and conclude that the arrangement of the two perpendicular lines that results in the maximum number of crossings for \(K_{m,n}^{\mathbf{DOL}}\) is when we have one of the numbers \(l, r\) equal to zero and one of the numbers \(t, b\) equal to zero, simultaneously (see Figure~\ref{fig_DOLmax}). Therefore, its number of crossings is as follows: 

\begin{equation*}
    cr\Big(K_{m,n}^{\mathbf{DOL}}\Big)_{\text{max}} = cr\Big(K_{m,n}^{((0,m),(n,0))}\Big) = \binom{m}{2} \cdot \binom{n}{2} = cr\Big(K_{m,n}^{\mathbf{DPL}}\Big).
\end{equation*}

In order to get the minimum number of crossings for \(K_{m,n}^{\mathbf{DOL}}\), we use inequality~\eqref{binom-min} in equation~\eqref{dol-compute-last}, i.e. we place the lines such that vertices on every line gets divided into two almost equal partitions. Hence, any of the following cases results in the minimum number of crossings for \(K_{m,n}^{\mathbf{DOL}}\):

\begin{equation*}
    \begin{cases}
      \lvert l-r\rvert \le 1,\\
      \text{}\\
      \lvert t-b\rvert \le 1.
    \end{cases}
\end{equation*}

Therefore, the formula for the minimum number of crossings for the \(K_{m,n}^{\mathbf{DOL}}\) is as follows:

\begin{equation}\label{dol-min}
    cr\Big(K_{m,n}^{\mathbf{DOL}}\Big)_{\text{min}} = \Bigg(\binom{\lfloor \frac{m}{2}\rfloor}{2} + \binom{\lceil \frac{m}{2}\rceil}{2}\Bigg) \cdot 
    \Bigg(\binom{\lceil \frac{n}{2}\rceil}{2} + \binom{\lfloor \frac{n}{2}\rfloor}{2}\Bigg). 
\end{equation}

By considering the parity of the positive integer \(k\), one can get that the following identity holds:

\begin{equation}\label{parity}
    \binom{\lfloor \frac{k}{2}\rfloor}{2} + \binom{\lceil \frac{k}{2}\rceil}{2} = \bigg\lfloor \frac{k}{2} \bigg\rfloor \cdot \bigg\lfloor \frac{k-1}{2} \bigg\rfloor. 
\end{equation}

Finally, using the expression~\eqref{parity} in equation~\eqref{dol-min} we get the following equation for the crossing number of \(K_{m,n}^{\mathbf{DOL}}\):

\begin{equation*}
    cr\Big(K_{m,n}^{\mathbf{DOL}}\Big)_{\text{min}} = \bigg\lfloor \frac{m}{2} \bigg\rfloor \cdot \bigg\lfloor \frac{m-1}{2} \bigg\rfloor \cdot \bigg\lfloor \frac{n}{2} \bigg\rfloor \cdot \bigg\lfloor \frac{n-1}{2} \bigg\rfloor.
\end{equation*}

\begin{figure}[H]
    \centering
    \includegraphics[width=6cm]{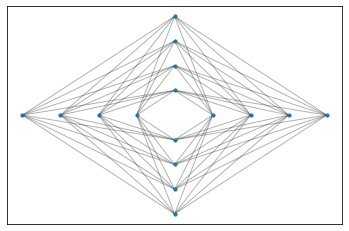}
    \caption{Double orthogonal lines arrangement of \(K_{8,8}:K_{8,8}^{\textit{DOL}}\) with minimum number of crossings}
    \label{fig_DOLmin}
\end{figure}

Figure~\ref{fig_DOLmin} shows the situation in which \(K_{m,n}^{\textit{DOL}}\) has its minimum number of crossings that is the case when 

\begin{equation*}
    K_{m,n}^{((\lfloor \frac{m}{2} \rfloor, \lceil \frac{m}{2} \rceil),(\lfloor \frac{n}{2} \rfloor, \lceil \frac{n}{2} \rceil))}.    
\end{equation*}

Zarankiewicz \cite{zarankiewicz1955problem} in its well-known conjecture states that the crossing number of the complete bipartite graph \(K_{m,n}\) is

\begin{equation*}
    \bigg\lfloor \frac{n}{2} \bigg\rfloor \cdot \bigg\lfloor \frac{n-1}{2} \bigg\rfloor \cdot \bigg\lfloor \frac{m}{2} \bigg\rfloor \cdot \bigg\lfloor \frac{m-1}{2} \bigg\rfloor.
\end{equation*}

Recall that we are considering the rectilinear crossing number, therefore, our formula stands for the rectilinear crossing number of \(K_{m,n}^{\textit{DOL}}\). If we assume that the Zarankiewicz conjecture is true, then the crossing number and rectilinear crossing number of the complete bipartite graph are equal; Kainen \cite{kainen1968problem}.

\subsection{Forth arrangement: Line Inside Circle LIC}\label{LIC-arrangement}
In this arrangement, we put part \(A\) of the complete bipartite graph \(K_{m,n}\) on an imaginary line inside the circle that contains the other part (i.e. part \(B\)) of the complete bipartite graph equally distanced on its curve, e.g., a regular polygon. We denote this arrangement by \(K_{m,n}^{\mathbf{LIC}}\). 

Suppose that the line (containing part \(A\) of the complete bipartite graph with \(m\) vertices) divides the circle into two partitions with \(k\) and \(l\) vertices (\(k+l=n\)). Then, we can count the number of crossings in the following way:

for every pair of vertices belonging to the part on the line, we get

\begin{equation}\label{one-pair-tpl-}
    \sum_{i=1}^{k-1} i + \sum_{j=1}^{l-1} j = \binom{k}{2} + \binom{l}{2}.
\end{equation}

And finally, since there exist \(\binom{m}{2}\) distinct pairs of vertices in the part on the line, we get the number of crossings as the following:

\begin{equation}\label{number-cr-tpl-}
    \binom{m}{2} \cdot \Bigg(\binom{k}{2} + \binom{l}{2}\Bigg).
\end{equation}

Since the function \(f(x) = x^2 + (z-x)^2\) attains its minimum at the point \(x=\lceil \frac{z}{2}\rceil\) or \(x=\lfloor \frac{z}{2}\rfloor\), the following inequality is true for all positive integers \(x,y\), where \(x+y=z\):

\begin{equation}\label{binom-min-}
    \binom{x}{2} + \binom{y}{2} \ge \binom{\lfloor \frac{z}{2}\rfloor}{2} + \binom{\lceil \frac{z}{2}\rceil}{2}.
\end{equation}

Moreover, for every two positive integers \(x,y\), where \(x \le y\), the following inequality is true:

\begin{equation}\label{m-n-inequality-}
    \binom{x}{2} \cdot \bigg \lfloor{\frac{y}{2}}\bigg \rfloor \cdot \bigg \lfloor{\frac{y-1}{2}}\bigg \rfloor \ge \binom{y}{2} \cdot \bigg \lfloor{\frac{x}{2}}\bigg \rfloor \cdot \bigg \lfloor{\frac{x-1}{2}}\bigg \rfloor.     
\end{equation}

\begin{figure}[H]
    \centering
    \includegraphics[width=6cm]{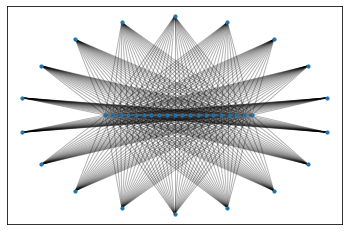}
    \caption{Line inside circle arrangement of \(K_{20,18}:K_{20,18}^{\mathbf{LIC}}\) with minimum number of crossings}
    \label{fig_LICmin}
\end{figure}

Hence, in order to get the minimum number of crossings in this arrangement, without loss of generality assume that \(m \ge n\). By using inequality~\eqref{m-n-inequality-} and inequality~\eqref{binom-min-} in equation~\eqref{number-cr-tpl-}, we get that we should put the part with larger number of vertices on the line (\(m\) vertices) and place the line on a diameter of the circle such that it divides the vertices on the circle into two parts with \(\lfloor \frac{n}{2} \rfloor\) and \(\lceil \frac{n}{2} \rceil\) vertices. Therefore, we get the following crossing number for \(K_{m,n}^{\mathbb{LIC}}\) arrangement (an example of such a case is shown in Figure~\ref{fig_LICmin}):

\begin{equation*}
    cr\Big(K_{m,n}^{\mathbf{LIC}}\Big) = \binom{m}{2} \cdot \bigg \lfloor{\frac{n}{2}}\bigg \rfloor \cdot \bigg \lfloor{\frac{n-1}{2}}\bigg \rfloor 
\end{equation*}

We believe that the arbitrary placement of vertices inside a circle, which we denote by \(K_{m,n}^{\mathbf{FIC}}\), does not bring anything new, so the minimum and maximum number of crossings would still be obtained in the cases when all the vertices inside the circle are located on a line. 

\begin{conjecture}
Let \(m \ge n\) and \(K_{m,n}^{\mathbf{FIC}}\) be an arrangement where \(m\) vertices of \(K_{m,n}\) are randomly distributed inside the circle that has \(n\) vertices on its curve. Then we have the following bound for its number of crossings:

\begin{equation*}\label{Cr-fic}
    \binom{m}{2} \cdot \bigg \lfloor{\frac{n}{2}}\bigg \rfloor \cdot \bigg \lfloor{\frac{n-1}{2}}\bigg \rfloor \le cr\Big(K_{m,n}^{\mathbf{FIC}}\Big) \le \binom{m}{2} \cdot \binom{n}{2} 
\end{equation*}
\end{conjecture}

\subsection{Double Circular arrangement: DC}\label{DC-arrangement} 

In this arrangement, we uniformly distribute the vertices of \(K_{m,n}\) on two concentric circles, i.e. one part on the \textit{inner} circle and the other part on the \textit{outer} circle (see Figure~\ref{double circular graphs}). 

Now we shall prove Proposition~\ref{Prop1}, but before starting the proof we give the following definitions.

\begin{definition}
We say a partitioning of the set \(V\) into \(V_1\) and \(V_2\) is a \emph{balanced partition} if \(\lvert \lvert V_1\rvert - \lvert V_2\rvert \rvert \le 1\).
\end{definition}

\begin{definition}
A line is called a \emph{balanced line} for the set \(V\) of vertices if it partitions the set into a balanced partition.
\end{definition}

\begin{proof}(Proposition~\ref{Prop1})
  
Let \(K_{q,s} = (X=U\cup V, E)\), where \(U=\{u_1, \cdots , u_q\}, V=\{v_1,\cdots , v_s\}\) and \(E=U \times V\). 
  
Let \(l_{i,j}\) be the line passing through vertices \(i,j \in U, i \not=j\). Line \(l_{i,j}\) partitions the set \(V\) of vertices on the outer circle into two parts, i.e. \(t_{i,j}\) and \(b_{i,j}\).
  
By using equation~\eqref{one-pair-tpl-} for every pair \(\{i,j\} \in U, i \not=j\), we compute the number of crossings by the following formula:
\begin{equation}\label{cross_compute}
    \sum_{1\le i < j \le q} \binom{\mid t_{i,j}\mid}{2} + \binom{\mid b_{i,j}\mid}{2},
\end{equation}
where \(\lvert t_{i,j}\rvert + \lvert b_{i,j}\rvert = q\), for all pairs \(\{i,j\} \in U, i \not=j\).
  
By the definition of the crossing number and using formula~\eqref{cross_compute}, we get the following equation:
  
  \begin{equation}\label{cross_def}
      cr\Big(K_{q,s}^{\mathbf{DC}}\Big) =  \min \Bigg( \sum_{1\le i < j \le q} \binom{\mid t_{i,j}\mid }{2} + \binom{\mid b_{i,j}\mid}{2} \Bigg),
  \end{equation}
   
where the minimum is taken over all sets of lines \(l_{i,j}\), where \(\{i,j\} \in U, i \not=j\).
  
It is trivial that the following inequality is always true:
  
  \begin{equation}\label{min_sum}
      cr\Big(K_{q,s}^{\mathbf{DC}}\Big) \ge \sum_{1\le i < j \le q} \min \Bigg( \binom{\mid t_{i,j}\mid}{2} + \binom{\mid b_{i,j}\mid}{2}  \Bigg) = A.
  \end{equation} 
  
Equality in \eqref{min_sum} occurs when all terms in \(A\) reach their minima simultaneously. 
    
Moreover, we know that every term inside the summation in equation~\eqref{min_sum} will be minimum only when the following holds:
  
\begin{equation}\label{min-condition}
    \lvert \mid t_{i,j}\mid - \mid b_{i,j}\mid \rvert \le 1 
\end{equation}

This means that every line \(l_{i,j}\) is a balanced line, i.e. equation~\eqref{min-condition} holds.  
  
In the hardest case, i.e. \(q,s \rightarrow \infty \), we just need that any line touching the inner circle to be a balanced line for the set of the vertices on the outer circle. Therefore, we fix the position of the outer circle and decrease the radius of the inner circle. Note that if any line \(l_{i,j}\) passes through a vertex \(p \in V\) that is on the outer circle, then we just rotate the inner circle a bit with respect to its center to avoid this to happen. By increasing the radius of the inner circle and rotating when needed, eventually at some point (i.e. sufficiently small radius for the inner circle) we reach a configuration where \textit{any line touching the inner circle will become a balanced line for the set of the vertices on the outer circle}. Since the \(s\) vertices on the outer circle form a regular \(s\)-gon and are equally distributed and distanced on its curve, they will be partitioned into a balanced partition as well, i.e. equi-cardinality division. 
  
Hence, we get the following equality: 
  
  \begin{equation}\label{cross_min}
       \sum_{1\le i < j \le q} \min \Bigg( \binom{\mid t_{i,j}\mid}{2} + \binom{\mid b_{i,j}\mid}{2} \Bigg) = \sum_{1\le i < j \le q} \binom{\lfloor \frac{s}{2}\rfloor}{2} + \binom{\lceil \frac{s}{2}\rceil}{2}.
  \end{equation}
  
  Equations~\eqref{cross_def}, \eqref{min_sum} and \eqref{cross_min} imply the following:
    
  \begin{equation*}
     cr\Big(K_{q,s}^{\mathbf{DC}}\Big) = \binom{q}{2} \bigg\lfloor\frac{s}{2} \bigg\rfloor \bigg\lfloor \frac{s-1}{2} \bigg\rfloor. 
  \end{equation*}

  Now consider the complete bipartite graph \(K_{m,n}\), and let \(m \ge n\), then by inequality~\eqref{m-n-inequality-} the following holds:
  
  \begin{equation*}
      \binom{m}{2} \bigg\lfloor\frac{n}{2} \bigg\rfloor \bigg\lfloor \frac{n-1}{2} \bigg\rfloor \le \binom{n}{2} \bigg\lfloor\frac{m}{2} \bigg\rfloor \bigg\lfloor \frac{m-1}{2} \bigg\rfloor.
  \end{equation*}

  Therefore, the number of crossings for \(K_{m,n}^\mathbf{DC}\) will be minimum if we put the part with \(m\) vertices on the inner circle and the part with \(n\) vertices on the outer circle, and this gives the following equality:

\begin{equation}\label{DC-min}
    cr\Big(K_{m,n}^\mathbf{DC}\Big) = \binom{m}{2} \cdot \bigg \lfloor{\frac{n}{2}}\bigg \rfloor \cdot \bigg \lfloor{\frac{n-1}{2}}\bigg \rfloor.
\end{equation}
\end{proof}

Now that we have the rectilinear crossing number of \(K_{m,n}^\mathbf{DC}\), a natural question to ask could be the following:

\textit{Let \(R\) and \(r\) be the radii of the outer and the inner circle, respectively. Does there exist a threshold for ratio \(\frac{R}{r}\) until which the number of crossings in double circular arrangement remains at its minimum?}

In the next section, we will answer the above mentioned question and introduce an appropriate problem setting to obtain that threshold.

\section{A geometric optimization problem}

We already know the rectilinear crossing number of \(K_{m,n}^\mathbf{DC}\) and that we can reach it if we choose sufficiently small radius for the inner circle. 

In this section, we will seek an answer to whether or not there exist a threshold for ratio \(\frac{R}{r}\) until which the equation~\eqref{DC-min} holds, where \(R\) is the radius of the outer circle and \(r\) is the radius for the inner one. Moreover, if such a minimum ratio exists, then we try to introduce a problem setting by which it is possible to obtain it. 

\subsection{Threshold: Maximum radius}\label{Threshold}

Let \(m\) and \(n\) be the number of vertices on inner and outer circles, respectively, where \(r\) and \(R\) are their radii and \(m \ge n\).

We define \(\theta\) as the angle between two radii of the same circle that are connected to two consecutive vertices on the same circle, then we will have:
\begin{equation}\label{thetas}
  \begin{cases}
    \theta_I = \frac{2\cdot \pi}{m}, &\quad \textit{for inner circle},\\
    \text{}\\
    \theta_O = \frac{2\cdot \pi}{n}, &\quad \textit{for outer circle}.
  \end{cases}    
\end{equation}

\begin{figure}[H]
    \centering
    \includegraphics[width=6cm]{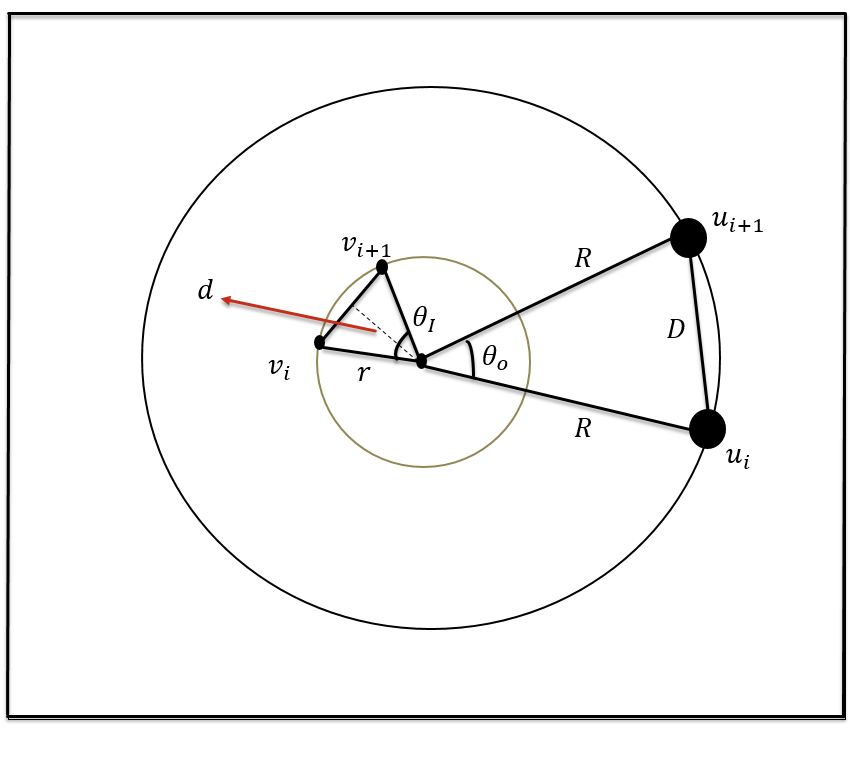}
    \caption{Parameters of two circles for \(K_{m,n}^\mathbf{DC}\)}
    \label{fig_parameters}
\end{figure}

We denote the length of the chord connecting two consecutive vertices on the outer circle by \(D = 2 \cdot R \cdot \sin{\big(\frac{\theta_O}{2}}\big)\). Figure~\ref{fig_parameters} shows the introduced parameters.

Let \(l_{i,j}\) be the line passing through pair \(\{i,j\}\) of vertices on the inner circle. And let line \(L_{i,j}\) be a line parallel and very close to \(l_{i,j}\) which lies between the center of the two circles and \(l_{i,j}\).   

\begin{figure}[H]
    \centering
    \includegraphics[width=5cm]{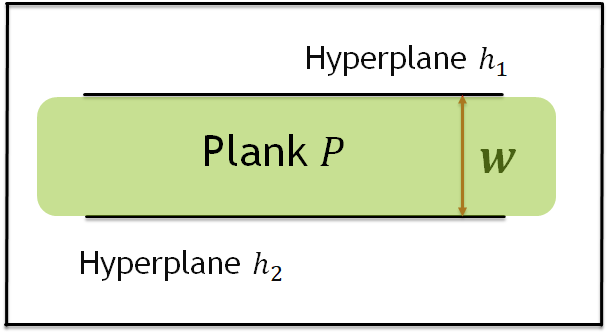}
    \caption{A plank \(P\) between lines \(h_1\) and \(h_2\) in \(\mathbb{R}^2\)}
    \label{plank}
\end{figure}

\begin{definition}

A space between two parallel hyperplanes \(h_1\) and \(h_2\) is called a \textit{plank}, where the distance between two hyperplanes is the width of plank and is denoted by \(w\). 
\end{definition}
 A plank \(P\) in \(2\)-dimensional space is shown in Figure~\ref{plank}.

Let plank \(P_i\) be defined by the line \(L_{i,i+ \lceil \frac{n}{2}\rceil -1}\) and its parallel line \(L_{i-1,i+\lceil \frac{n}{2}\rceil}\) where indices of the lines stand for vertices on the outer circle. Note that every one of these lines is sufficiently close to the vertices by which they are defined, but does not pass through the vertices. Then, we have the following definition for a \textit{main plank} that is specific for this work.

\begin{figure}[H]
    \centering
    \includegraphics[width=5cm]{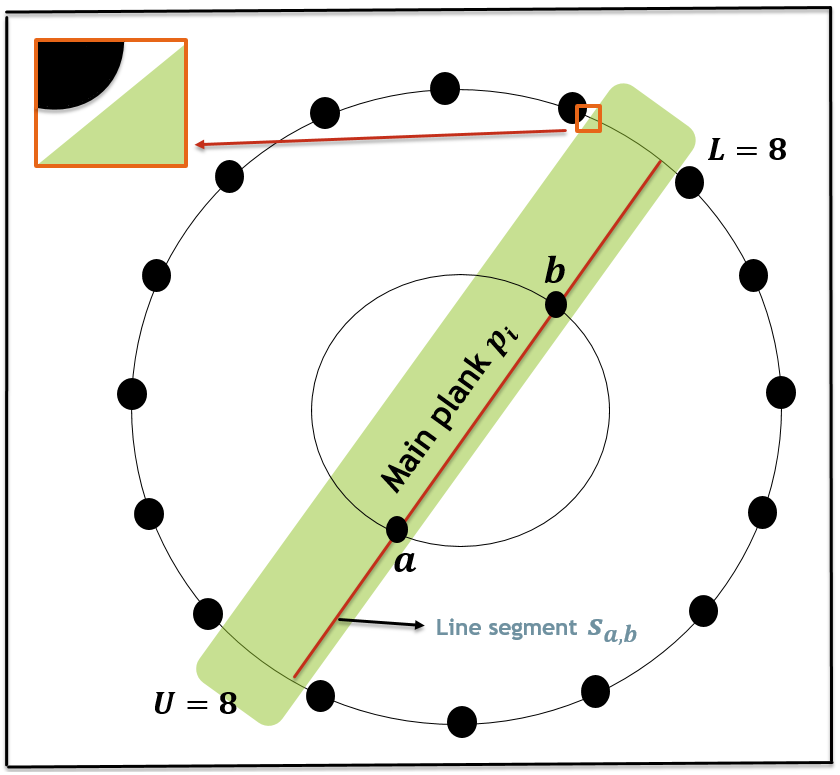}
    \caption{A main plank \(p_i\) and a line segment \(s_{a,b}\) which belongs to it} 
    \label{main-plank}
\end{figure}

\begin{definition}
A \textit{main plank} \(p_i\) is a plank that contains the center of the circles and its defining lines are parallel to the diameter of the circles. Moreover, it has two properties; first, it divides the vertices on the outer circle into two almost equal parts and second, it does not contain any of the vertices on the outer circle (see Figure~\ref{main-plank}).
\end{definition}

Denote the set of all main planks by \(M = \{p_1,\cdots,p_{\lceil \frac{n}{2}\rceil -1 }\}\), where their widths \(w_i\) are all the same and \(w_i + \varepsilon = D = 2 \cdot R \cdot \sin{\big(\frac{\theta_O}{2}}\big)\), where \(\varepsilon > 0\) and is sufficiently small.

Denote by \(s_{a,b}\) the line segment with its end points being the intersection of the line \(l_{a,b}\) with the outer circle. Then we have the following definition. 

\begin{definition}
A line segment belongs to a plank if it is totally inside the plank, i.e. it lies inside the plank and its intersection with the plank's lines is either a subset of its end points or the empty set. 
\end{definition}

To achieve the minimum number of crossings, we have the following condition: every line segment \(s_{a,b}\) that corresponds to the pair \(\{a,b\}\) of the vertices on the inner circle must belong to at least one main plank. This is equivalent to, each line \(l_{a,b}\) is a balanced line for the set of the vertices on the outer circle (i.e. divides the vertices on the outer circle into two almost equal parts). 

Now we can define the following \textit{geometric optimization problem}: for given \(R, m, n\) and rotation of the outer (or inner) circle, determine maximum value of \(r\) for which the number of crossings of \(K_{m,n}\) in double circular arrangement is minimized. 

Let \(m \ge n\), then by Proposition~\ref{Prop1} we know that to get the minimum number of crossings in double circular arrangement of \(K_{m,n}\), we must put the parts with \(m\) and \(n\) vertices on the inner and outer circle, respectively. Denote the two parts of the vertices by \(V\) and \(U\), where \(V=\{v_1,\cdots , v_m\}\) and \(U=\{u_1, \cdots , u_n\}\). Let line segment \(s_{a,b}\) be the intersection of the line \(l_{a,b}\), that passes through vertices \(a\) and \(b\) belonging to the set \(V\), with the outer circle and \(r\) be the radius of the inner circle. Then the following is the statement of the geometric optimization problem: 

\begin{maxi!}|l|[3]
{}{r}{}{}
\addConstraint{\forall s_{a,b}, \exists p_i, \textit{such that}: s_{a,b} \subset p_i, &\quad \forall a,b \in V, a \not= b, p_i \in M}
\addConstraint{s_{a,b} \cap U = \varnothing
, & \quad \forall a,b \in V, a \not= b}{}
\addConstraint{r > 0}{}
\end{maxi!}

Recall that the solution to this optimization problem is a threshold for the radius of the inner circle such that the number of crossings of \(K_{m,n}^\mathbf{DC}\) is minimum, we call it \textit{crossing threshold}. Now that we have the crossing threshold, the following question raises:

\textit{How the number of crossings of \(K_{m,n}^\mathbf{DC}\) will change if we pass its crossing threshold?}

In the next section, we will prove Lemma~\ref{Lemma1} that provides a partial answer to this question.

\section{Above the Crossing Threshold}\label{above the threshold}

We call the maximum \(r\) after which the number of crossings stops to be minimum, the \textit{crossing threshold} and denote it by \(T_{\text{cr}}\).

Now we will prove the Lemma~\ref{Lemma1} which gives a bound for the number of crossings when \(r > T_{\text{cr}}\). Note that when \(r < T_{\text{cr}}\), Lemma~\ref{Lemma1} gives the exact number of crossings which is the same as the crossing number of \(K_{m,n}^\mathbf{DC}\) that we obtained in Proposition~\ref{Prop1}.

\begin{proof}(Lemma~\ref{Lemma1})\\

  In \(K_{m,n}^\mathbf{DC}\), let \(m\) and \(n\) vertices be on the inner and the outer circles, respectively, where \(r\) and \(R\) are their radii.

  Let line \(l_{a,b}\) be the line passing through pair \(\{a,b\}\) of vertices on the inner circle that divides the vertices on the outer circle into two partitions. Let \(U_{a,b}\) and \(L_{a,b}\) be the number of vertices that are in different sides of the line \(l_{a,b}\) (i.e. the number of the vertices in the partitions caused by line \(l_{a,b}\)), and let \(d_{a,b}\) be the distance between line \(l_{a,b}\) with the center of the circles (or with the diameter that is parallel to this line).

  We partition the set of all lines passing through pairs of vertices of the inner circle into different classes, where every class \(C_j\) is defined by an arbitrarily chosen member of it which we call its \textit{representative line} \(r_{j}\). Moreover, a line \(l_{a,b}\) belongs to the class \(C_j\) if and only if \(b-a=j\), where \(j \in \{1, \cdots, (\lceil \frac{m}{2}\rceil -1)\}\), some of the line classes are shown in Figure~\ref{fig_lc}. 
  
  Let \(d_j\) be the distance of the representative line \(r_j\) with the center of the circles. Then the following holds: 
  
  \begin{equation}\label{distance-dj}
      d_j = r \cdot \cos{\Bigg(\frac{\pi - (j\cdot \theta_I)}{2}\Bigg)} = r \cdot \sin{\Bigg(\frac{j\cdot \theta_I}{2}\Bigg)}, 
  \end{equation}

  where from equation~\eqref{thetas}, \(\theta_I = \frac{2\cdot \pi}{m}\).
  
  Line \(r_j\) divides the vertices on the outer circle into two partitions. Let \(U_j\) and \(L_{j}\) be the number of the vertices that lie in these two partitions caused by the line \(r_j\).
  
  By the symmetry of the \(K_{m,n}^\mathbf{DC}\), the following properties hold for every class \(C_j\) of lines:
  
  \begin{itemize}
      \item \(d_k=d_j, \forall k \in C_j\) (i.e. The distance \(d_k\) of all the lines \(l_k\) belonging to the same class \(C_j\) are the same and equal to \(d_j\)), 
      \item \((U_k,L_k)=(U_j,L_j), \forall k \in C_j\) (i.e. all pairs \(\{U_k,L_k\}\) corresponding to different lines \(l_k\) are the same and equal to \(\{U_j, L_j\}\)).
  \end{itemize}

  \begin{figure}[H]
    \centering
    \includegraphics[width=6cm]{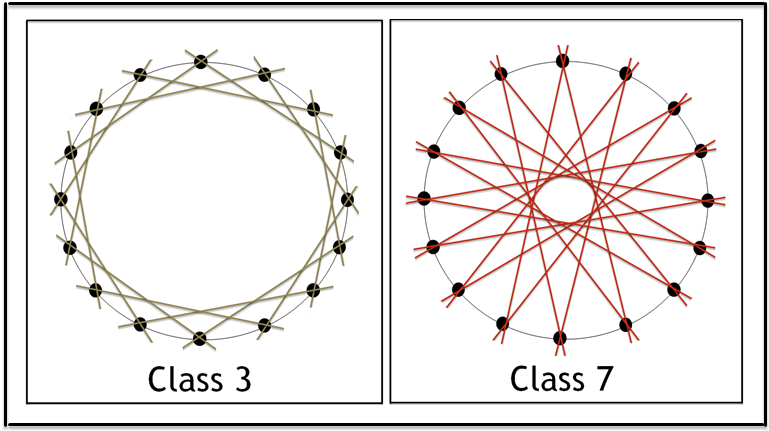}
    \caption{Line classes \(C_3\) and \(C_7\) for \(m=16\)}
    \label{fig_lc}
  \end{figure}

  Therefore, the number of the crossings corresponding to all the lines of the class \(C_j\) are equal. Hence, we can estimate the total changes in the number of crossings for \(K_{m,n}^\mathbf{DC}\) through estimating changes in the number of crossings corresponding to each representative line \(r_j\). Furthermore, every class of the lines contains exactly \(m\) lines. Therefore, the following equations hold:
  
  \begin{equation}
      \begin{cases}
        \dbinom{m}{2} - \bigg( m \cdot \Big(\lceil \frac{m}{2}\rceil -1\Big)\bigg) = {\frac{m}{2}}, &\quad \textit{for even \(m\)},\\
        \text{}\\
        \dbinom{m}{2} - \bigg( m \cdot \Big(\lceil \frac{m}{2}\rceil -1\Big)\bigg) = 0, &\quad \textit{for odd \(m\)}.
      \end{cases}
  \end{equation}

  This means, if \(m\) is even, then there exist some lines that are not included in any of the line classes. However, these \(\frac{m}{2}\) lines form the class \(C_{\frac{m}{2}}\) themselves. This class has a special property that is, every line in this class contains a main diameter of the inner circle. Therefore, each of them is a balanced line (i.e. divides the vertices on the outer circle into two almost equal parts). Then, we get the following number of crossings corresponding to these lines:
  
  \begin{equation}\label{term-even-m}
      \frac{m}{2} \cdot \bigg\lfloor\frac{n}{2} \bigg\rfloor \bigg\lfloor \frac{n-1}{2} \bigg\rfloor.
  \end{equation}

  And the number of crossings corresponding to the \((\lceil \frac{m}{2}\rceil -1 )\) classes can be counted as follows:

  \begin{equation}\label{common-term}
      m \cdot \sum_{j=1}^{\lceil \frac{m}{2}\rceil -1} \binom{U_j}{2} + \binom{L_j}{2} .
  \end{equation}
 
  Now we shall see how \(U_j\) and \(L_j\) will change by increasing the radius \(r\). Recall that when \(r < T_{\text{cr}}\), we have the following:
  
  \begin{equation}\label{ul-opt}
        \lvert U_j - L_j \rvert \le 1. 
  \end{equation}

Moreover, we know that every line segment \(s_{a,b}\) must belong to at least one main plank, that gives \(d_{a,b} < D\) for all pair \(\{a,b\}\) of vertices on the inner circle. The same is true for all the representative lines \(r_j\). 

Therefore, if we increase the radius \(r\), then the line \(r_j\) moves in a way that \(d_j\) increases, and then equation~\eqref{ul-opt} will not hold anymore and instead we will have the following: 
  
  \begin{equation}\label{after-threshold}
      \begin{cases}
        U_{j} = \frac{n}{2} - (c \cdot \beta_{j}), L_{j} = \frac{n}{2} + (c \cdot \beta_{j}), &\quad \textit{for even \(n\)},\\
        \text{}\\
        U_{j} = \frac{n-1}{2} - (c \cdot \beta_{j}), L_{j} = \frac{n+1}{2} + (c \cdot \beta_{j}), &\quad \textit{for odd \(n\)},
      \end{cases}
  \end{equation}
  
  where \(\beta_j\) satisfies the following inequality:
  
  \begin{equation*}\label{beta-j}
      \beta_j \cdot D < d_j < (\beta_j+1) \cdot D,
  \end{equation*}
  
  and \(c \in \{1,2\}\). In order to find \(c\), recall that every line segment belongs to at least one main plank, denote the main plank that contains line segment \(s_{a,b}\) by \(p_{i}^{a,b}\). Then, we can get \(c\) from the following cases (see Figure~\ref{c-case1}):
  
  \begin{equation*}\label{c-beta}
      \begin{cases}
        c = 1, &\quad \textit{if \(r_j\) is not parallel to \(p_{j,i}\)},\\
        \text{}\\
        c = 2, &\quad \textit{if \(r_j\) is parallel to \(p_{j,i}\)}.
      \end{cases}
  \end{equation*}
  
  
  \begin{figure}[H]
    \centering
    \includegraphics[width=8cm]{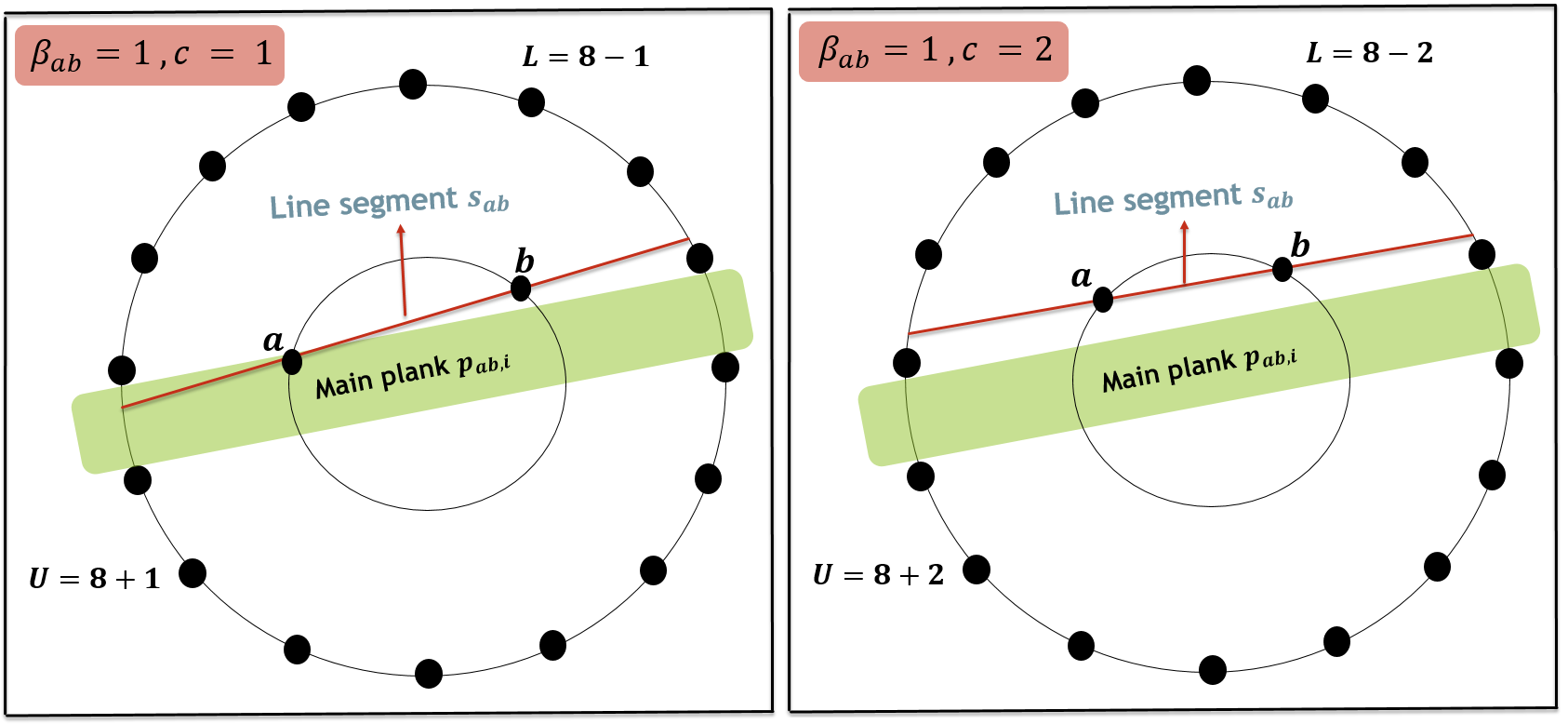}
    \caption{left: \(s_{a,b}\) not parallel to \(p_{ab,i}\) (\(c=1\)) , right: \(s_{a,b}\) parallel to \(p_{ab,i}\) (\(c=2\))}
    \label{c-case1}
  \end{figure}
  
  By using equation~\eqref{after-threshold} in equation~\eqref{common-term}, we get the following:
  
  \begin{equation*}
      \begin{cases}
        m \cdot \sum_{j=1}^{\lceil \frac{m}{2}\rceil -1} \Bigg (\dbinom{\frac{n}{2} - (c \cdot \beta_{j})}{2} + \dbinom{\frac{n}{2} + (c \cdot \beta_{j})}{2} \Bigg ), &\quad \textit{for even \(n\)},\\
        \text{}\\
        m \cdot \sum_{j=1}^{\lceil \frac{m}{2}\rceil -1} \Bigg (\dbinom{\frac{n-1}{2} - (c \cdot \beta_{j})}{2} + \dbinom{ \frac{n+1}{2} + (c \cdot \beta_{j})}{2} \Bigg ), &\quad \textit{for odd \(n\)},\\
      \end{cases}
  \end{equation*}

  And by doing some calculations we will have the following:
  
  \begin{equation}\label{final-term}
    \begin{cases}
     m \cdot \sum_{j=1}^{\lceil \frac{m}{2}\rceil -1} \bigg ( \Big(\frac{n}{2}\Big)^2 - \Big(\frac{n}{2}\Big) + (c\cdot \beta_j)^2 \bigg ), &\quad \textit{for even \(n\)},\\
     \text{}\\
     m \cdot \sum_{j=1}^{\lceil \frac{m}{2}\rceil -1} \bigg ( \Big( \frac{n-1}{2} \Big)^2 + (c\cdot \beta_j)^2 + (c\cdot \beta_j) \bigg ), &\quad \textit{for odd \(n\)},\\
    \end{cases}
  \end{equation}

  Now assume that none of the representative lines \(r_j\) (also we can make the assumption for all line segments) are parallel to their main planks \(p_{j,i}\), then we put \(c=1\) and we get lower bounds for all parity combinations of \(m\) and \(n\).
  
  Recall that for even \(m\), we have a term that we obtained in equation~\eqref{term-even-m}. By doing some simplifications on equations~\eqref{final-term}, we get the following lower bounds:
  
  \begin{equation*}
      N_{cr}\Big(K_{m,n}^\mathbf{DC} \Big) \ge
      \begin{cases}
        \dbinom{m}{2} \cdot {\bigg(\frac{n^2}{4}-\frac{n}{2}\bigg) + m\cdot \sum_{j=1}^{\frac{m}{2}-1} (\beta_j)^2}, &\quad \textit{for even \(n\)},\\
         \text{}\\
        \dbinom{m}{2}\cdot {\bigg(\frac{n-1}{2} \bigg)^2 + m\cdot \sum_{j=1}^{\frac{m-1}{2}} \beta_j + (\beta_j)^2} ,&\quad \textit{for odd \(n\)},\\
      \end{cases}
  \end{equation*}
  
  Now assume that all representative lines \(r_j\) (also we can make the assumption for all line segments) are parallel to their main planks \(p_{j,i}\), then we put \(c=2\) and we get upper bounds for all parity combinations of \(m\) and \(n\).
  
  Recall that for even \(m\), we have a term that we got in equation~\eqref{term-even-m}. By doing some simplifications on equations~\eqref{final-term}, we get the following upper bounds:

  \begin{equation*}
      N_{cr}\Big(K_{m,n}^\mathbf{DC} \Big) \le
      \begin{cases}
        \dbinom{m}{2} \cdot {\bigg(\frac{n^2}{4}-\frac{n}{2}\bigg) + 4m\cdot \sum_{j=1}^{\frac{m}{2}-1} (\beta_j)^2}, &\quad \textit{for even \(n\)},\\
         \text{}\\
        \dbinom{m}{2}\cdot {\bigg(\frac{n-1}{2} \bigg)^2 + 2m\cdot \sum_{j=1}^{\frac{m-1}{2}} \beta_j + 2\cdot(\beta_j)^2} ,&\quad \textit{for  odd \(n\)},\\
      \end{cases}
  \end{equation*}

  Now we can merge all these bounds and get the following general bounds:
  
\begin{equation*}
    \begin{cases}
        N_{cr}\Big(K_{m,n}^{\mathbf{DC}}\Big) - \binom{m}{2} \bigg\lfloor\frac{n}{2}\bigg\rfloor  \bigg\lfloor\frac{n-1}{2}\bigg\rfloor \ge m \sum_{j=1}^{\lceil\frac{m}{2}\rceil-1} \Bigg(\Big(\frac{1-(-1)^n}{2}\Big)\beta_j + (\beta_j)^2\Bigg)   
        \text{ }\\\\
        N_{cr}\Big(K_{m,n}^{\mathbf{DC}}\Big) - \binom{m}{2} \bigg\lfloor\frac{n}{2}\bigg\rfloor  \bigg\lfloor\frac{n-1}{2}\bigg\rfloor \le 4m  \sum_{j=1}^{\lceil\frac{m}{2}\rceil-1} \Bigg(\Big(\frac{1-(-1)^n}{4}\Big)\beta_j + (\beta_j)^2\Bigg).
    \end{cases}
\end{equation*}  
  
  This concludes the proof.

\end{proof}

\section{Conclusions}

We have obtained the minimum number of crossings for a specific drawing of the complete bipartite graph (i.e. two concentric circles that contain two parts of the \(K_{m,n}\)), stated in Proposition~\ref{Prop1}. We also introduced a geometrical optimization problem to obtain a rotation-specific threshold for the ratio of radii of the circles for the minimum number of crossings. Finally, we found bounds for changes in the number of crossings for ratios of the radii larger than the obtained threshold.

As a question for further investigation, we propose the following:

\textit{for given \(m,n,r\) and \(R\), what is the set of the optimal angles for rotation which minimizes the number of crossings?}

\section{Acknowledgments}

I thank Prof. Elena Bazanova for her great help and comments on the text of the manuscript.



\bibliography{./references}
\end{document}